# The Visual Pattern in the Collatz Conjecture and Proof of No Non-Trivial Cycles

By Fabian S. Reid | fsreid@wgu.edu

**Abstract** We present the long sought visual pattern in the Collatz problem with the aid of the logarithmic spiral $r = 2 \cdot 2^{\theta/\pi}$. Using this newly discovered pattern, we show that the Collatz problem is linked to primes via Jacobsthal numbers. We then prove that no non-trivial cycles exist on the spiral and by extension in the Collatz conjecture.

## 1. Introduction

The Collatz Problem (also called the Collatz conjecture or 3x + 1 Problem), named after Lothar Collatz who first proposed the conjecture in 1937[1], is notorious for its contrasting nature of being simple to understand yet difficult to prove and, as stated in [2], 'absorbing massive amounts of time from both professional and amateur mathematicians'. It states that for any natural number, if the number is even, divide it by 2, but if the number is odd, then multiply it by 3 and add 1; at some point during this binary process the answer will be 1. With no pattern to work with, this problem is often compared to the rising and falling movements of a hailstone, chaotic and unpredictable. Some quotes on the problem reveal just how elusive it is:

> "For about a month everyone at Yale worked on it, with no result. A similar phenomenon happened when I mentioned it at the University of Chicago." – Shizuo Kakutani, 1960. [2]
>
> "Mathematics is not yet ripe enough for such questions." – Paul Erdős, 1983. [3]
>
> "This is an extraordinarily difficult problem, completely out of reach of present day mathematics." – Jeff Lagarias, 2010. [4]

To date, all numbers below $10^{20}$ have been verified to converge to 1 with the help of a network of computers.[5]

## 2. Connecting the Dots to Arrive at the Visual Pattern

Let T: $\mathbb{N} \longrightarrow \mathbb{N}$ be the Collatz function, defined as

$$T(n) = \begin{cases} \frac{n}{2} & \text{if n is even,} \\ 3n+1 & \text{if n is odd.} \end{cases} \quad (1.1)$$

Let $T^i(n)$ be the output of the $i$th iteration of $T$ on $n$. The Collatz conjecture states that for all $n$, there exists some $i$ such that $T^i(n) = 1$. If $n$ is even, then under the Collatz function (1.1) we would repeatedly divide by 2 until the result is odd. If $n$ is an odd multiple of 3, that is, $n = 3(2x-1)$ for $x \epsilon \mathbb{N}$, then multiplying $3(2x-1)$by 3 and adding 1 gives $9(2x-1)+1 = 18x-8$. Since $18x-8$ is not a multiple of 3, then repeated division by 2 will at some point give an odd number that is not a multiple of 3. However, odd numbers that are not multiples of 3 are all of the form $6x \pm 1$ for $x \epsilon \mathbb{N}$. This leads to Lemma 1.1.

**Lemma 1.1** If the Collatz Conjecture is true for all natural numbers of the form $6x \pm 1$ where $x \epsilon \mathbb{N}$, then it is true for all natural numbers.

We will now delve into uncharted territory to discover the long sought visual pattern in the Collatz conjecture.

Consider the set of odd numbers in the form $2n - 1$. The first iterate under the Collatz function (1.1) is $6n - 2$. If we halve the first iterate of each odd number, starting from 1, we get the arithmetic sequence 2, 5, 8, 11, 14,… $3n - 1$ (A016789)[6]. If we divide each term of the sequence A016789 by $2^k$ for $k \in \mathbb{N}$ and k as large as possible to achieve an odd number, then we get the sequence 1, 5, 1, 11, 7, 17, 5, 23, 13, 29, 1, 35, 19, 41, 11, 47, 25, 53, 7, 59, 31, 65, 17, 71, 37, 77, 5, 83, 43, 89, 23, 95, 49, 101, 13, 107, 55, 113, 29, 119, 61, 125, 1, 131, 67, 137, 35…. (A075677). Each term in A075677 occurs infinitely many times in a spiraling pattern, with 1 starting a new revolution as described below. We will call each revolution on the spiral a level.

The first level contains the first two terms, 1 and 5. The second level contains the next 8 terms: 1, 11, 7, 17, 5, 23, 13, and 29. The third level contains the next 32 terms: **1**, 35, 19, 41, **11**, 47, 25, 53, **7**, 59, 31, 65, **17**, 71, 37, 77, **5**, 83, 43, 89, **23**, 95, 49, 101, **13**, 107, 55, 113, **29**, 119, 61, and 125. Note that all the elements of the first level are present in the second level and demarcate the second level into distinct groups of three new elements. Likewise, all the elements of the second level are present in the third level (in bold print) and demarcate the third level into distinct groups of three new elements. We note that the number of elements at each level is an odd power of 2; the first level has $2^1$ elements, the second level has $2^3$ elements, the third level has $2^5$ elements, and so on. The levels continue in this pattern indefinitely.

Consider the sequence A075677. If we let $a(n) = \frac{1 - A075677(n)}{6}$ if 6|(A075677(n)-1) and $a(n) = \frac{A075677(n) + 1}{6}$ if 6|(A075677(n)+1), then we get the sequence 0, 1, 0, 2, -1, 3, 1, 4, -2, 5, 0, 6, -3, 7, 2, 8, -4, 9, -1, 10, -5, 11, 3, 12, -6, 13, 1, 14, -7, 15, 4, 16, -8, 17, -2, 18, -9, 19, 5, 20, -10, 21, 0, 22, -11, 23,…( A329480), which was first discovered by the author.

Next, note that the sequence A016789 can be generated by the polar logarithmic function $r(\theta_n) = 2 \cdot 2^{\theta_n/\pi}$ where $\theta_n = \ln\left(\frac{3n-1}{2}\right) \cdot \frac{\pi}{\ln(2)}$ for $n \in \mathbb{N}$.
For example, $\theta_1 = \ln\left(\frac{3\cdot 1-1}{2}\right) \cdot \frac{\pi}{\ln(2)} = \ln(1) \cdot \frac{\pi}{\ln(2)} = 0$ which means $r(\theta_1) = r(0) = 2 \cdot 2^{0/\pi} = 2$, the first term of the sequence. Likewise, $\theta_2 = \ln\left(\frac{3\cdot 2-1}{2}\right) \cdot \frac{\pi}{\ln(2)} = \ln\left(\frac{5}{2}\right) \cdot \frac{\pi}{\ln(2)} = 4.15296$ to five decimal places which means $r(\theta_2) = r(4.15296) = 2 \cdot 2^{4.15296/\pi} = 5$, the second term of the sequence. Therefore, if we graph $r_n = 2 \cdot 2^{\theta_n/\pi}$ and add rays at angles $\theta_n$, such that $\theta_n = \ln\left(\frac{3n-1}{2}\right) \cdot \frac{1}{\ln(2)}$ for $n \epsilon \mathbb{N}$ and extend each ray from $r_n$ to infinity, we get Figure 1.1 below. The polar grid is removed for simplicity. Note that between any two successive rays at one level, there are three rays at the immediate next level; this reflects the previously stated observation that elements of the sequence A075677 at one level demarcate the next level into distinct groups of three elements at the next level.

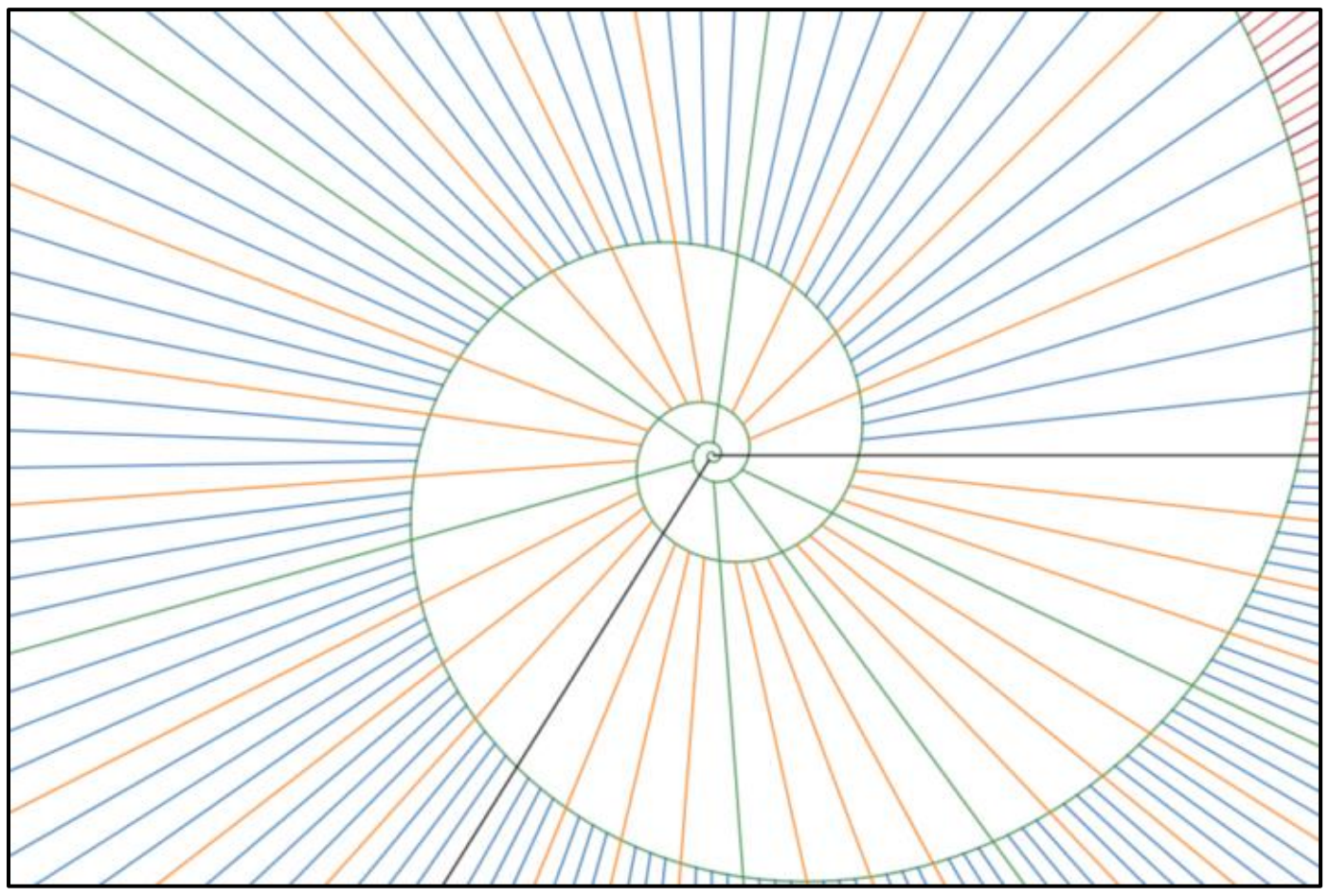

Figure 1.1. The polar graph of $r_n$ with extended rays along $\theta_n$.

If we label each ray, $\theta_n$, with the corresponding term from A329480 we get Figure 1.2 below, albeit in a smaller window.

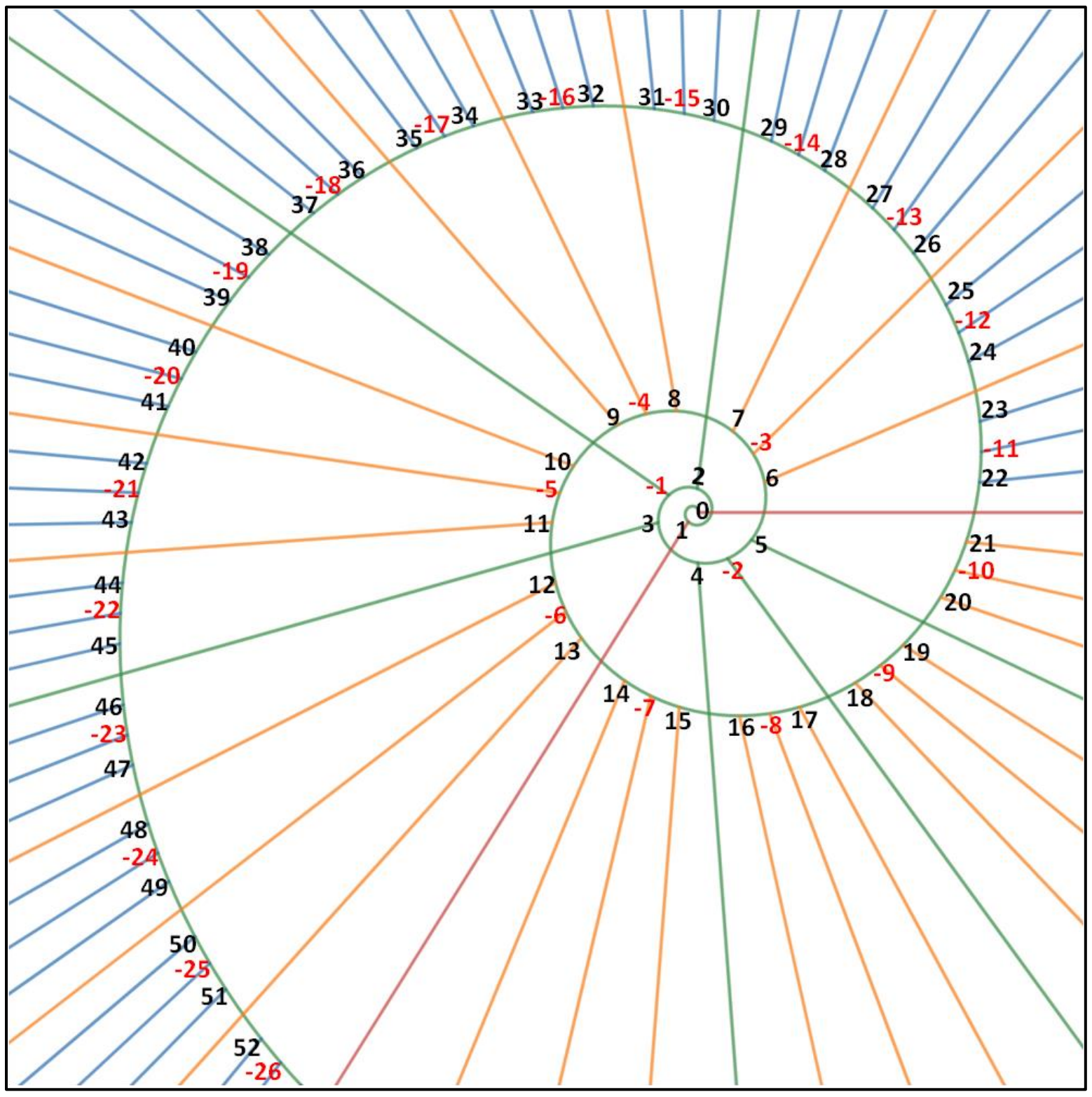

Figure 1.2. Figure 1.1 with each ray, $\theta_n$, labeled with the corresponding term from A329480.

**Remark 1.1** A level in Figure 1.2 corresponds to an angle of $2\pi$. The $kth$ level of Figure 1.2 contains integers of A329480 that fall in the interval $2\pi(k-1) \leq \theta < 2\pi k$. A number falls $k$ levels if it is moved backward through an angle of $2\pi k$ for $k \in \mathbb{N}$.

The first level of Figure 1.2 has only 0 and 1. The second level has 2, -1, 3, 4, -2 and 5.

Note that each positive number, $x$, on the spiral in Figure 1.2 represents the natural number $6x - 1$ and each negative number, $x$, represents the natural number $6(-x) + 1$. Since 65 = 6·11−1 is of the form 6x−1, it is represented by 11 on the spiral in Figure 1.2. Since 61 = 6·10+1 is of the form 6x+1, it is represented by −10 on the spiral in Figure 1.2.

**Proposition 1.1** The spiral trajectory of $6x \pm 1$ can be determined by Figure 1.2 by moving forward $x$ rays if $x$ is positive, moving backward $x$ rays if $x$ is negative, and continuing this forward-backward movement until '0' is reached.

Examples. Since $65 = 6 \cdot 11 - 1$, starting at 11 we move 11 rays forward to arrive at -8. Moving 8 rays backward we arrive at -6. Moving 6 rays backward we arrive at -1. Moving 1 ray backward from -1 we arrive at 2. Continuing this forward-backward movement we get the Spiral trajectory of 65 to be 11 → -8 → -6 → -1 → 2 → 3 → -2 → 1 → 0. If we convert each integer in the spiral trajectory of 65 to its corresponding number in the form $6x \pm 1$ we get 65 → 49 → 37 → 7 → 11 → 17 → 13 → 5 → 1, which are all the odd numbers in the Collatz trajectory of 65. Similarly, the spiral trajectory of 61 is -10 → 4 → 6 → 9 → 1 → 0. The corresponding numbers in the form $6x \pm 1$ give 61→ 23 → 35 → 53 → 5 → 1, which are the odd numbers in the Collatz trajectory of 61.

**Theorem 1.1** Let n be a natural number of the form $6x \pm 1$. Then the numbers in the spiral trajectory of n converted to their corresponding numbers in the form 6x±1 yield the odd numbers in the Collatz trajectory of $n$.

**Proof :** The sequence (A329480) and Figure 1.2 are developed from the sequence 1, 5, 1, 11, 7, 17, 5, … (A075677) which in turn is developed from the sequence 2, 5, 8, 11, 14, 17, 20, 23,… (A016789). The sequence A016789 can be expressed explicitly as $a(n) = 3n - 1$. Numbers of the form $6x - 1$ are the even indexed terms of A016789. Expressing $6x - 1$ in the form $3n - 1$ gives

$$3(2x) - 1. \tag{1.2}$$

Under the Collatz function $6x - 1$ becomes 3(6x – 1) +1 = 18x – 2. Dividing by 2 based on (1.1) yields $9x - 1$ which in the form $3n - 1$ equals

$$3(3x) - 1. \tag{1.3}$$

Since $n = 2x$ in (1.2) and $n = 3x$ in (1.3) we see that under the Collatz function, $6x - 1$ is moved x terms forward in A016789. Therefore, the numbers of the form $6x - 1$ are represented by positive integers in Figure 1.2 and correspond to moving forward $x$ rays on the spiral.

Now the $4i + 1$ indexed terms in A016789 can be divided by 2 to get the corresponding terms in A075677 which are all of the form $6x + 1$. This means that the $4i + 1$ indexed terms in A075677 equal to $12x + 2$ in A016789. If we write $12x + 2$ in the form $3n - 1$ we get

$$\begin{aligned} & 12x + 2 \\ = & 12x + 3 - 1 \\ = & 3(4x + 1) - 1 \end{aligned} \tag{1.4}$$

Under the Collatz function $6x + 1$ becomes $3(6x + 1) + 1 = 18x + 4$. Dividing by 2 gives $9x + 2$ which when written in the form $3n - 1$ gives

$$\begin{aligned} & 9x + 2 \\ = & 9x + 3 - 1 \\ = & 3(3x + 1) - 1. \end{aligned} \tag{1.5}$$

Since $n = 4x + 1$ in (1.4) and $n = 3x + 1$ in (1.5), we see that under the Collatz function, $6x + 1$ is moved $x$ terms backward in the sequence A016789. Therefore, the numbers of the form $6x + 1$ are represented by negative integers in Figure 1.2 and correspond to moving $x$ rays backward on the spiral. ∎

The conclusion of all this is that the spiral in Figure 1.2 is a model illustrating the visual pattern in the Collatz trajectory of numbers of the form $6x \pm 1$. This visual pattern is in stark contrast to the best-known visuals for the Collatz conjecture, as can be seen on Wikipedia, which show no pattern and are often described as chaotic and unpredictable [7].

## 3. Collatz Conjecture and Primes

The forward and backward movements in Figure 1.2 can be summarized as follows:

1. If we start at a **positive even number** we will continue moving forward until we get to an odd number. This odd number will always be a multiple of 3. That is, if $x$ is a positive even number, $x$ will map unto $\frac{3}{2}x$. This mapping loops until we arrive at an odd multiple of 3. For example, $4 \to 6 \to 9$ and $8 \to 12 \to 18 \to 27$.
2. If $x$ is a **positive odd number** we have two possibilities for an output: either x will map unto a negative integer iff $x \equiv 3 \mod 4$ or $x$ will fall 1 or more levels to either a positive or negative value which is governed by an infinite set of Diophantine equations, which we will discuss later.
3. If $x$ is a **negative odd number**, then $x$ will map unto a positive number given by $\frac{-3}{2}x + \frac{1}{2}$.
4. If $x$ is a **negative even number,** as is the case for a positive odd number, we have infinite possibilities for the output.

In a Collatz trajectory, multiplying an odd number by 3 and adding 1 will always be followed by division of 2. However, multiplying an odd number by 3 and adding 1 may also necessitate two consecutive divisions by 2 (or division by 4), or three consecutive divisions by 2, (or division by 8), and so on. Herein lies the infinite set of Diophantine equations that govern the movements of odd positive integers and even negative integers in Figure 1.2. We need a different equation for each $2^k$ division that we might need to make.

We now summarize the forward or backward movement along the spiral in Figure 1.2 in the following Lemma.

**Lemma 1.2** Let $x$ be a number on the spiral. Then there exists $k \in \mathbb{N}$ such that the resulting integer output of the forward (if $x > 0$) or backward (if $x < 0$) movement is given by the Diophantine equation

$$f_k(x) = \frac{(-1)^{k-1} \cdot 3}{2^k} x + \frac{-2 \cdot (-1)^{k-1} + 2^k}{6 \cdot 2^k} \qquad \text{if } x > 0 \tag{1.6}$$

or

$$f_{-k}(x) = \frac{(-1)^{k} \cdot 3}{2^k} x + \frac{4 \cdot (-1)^{k-1} + 2^k}{6 \cdot 2^k} \qquad \text{if } x < 0 \tag{1.7}$$

Examples

1. Starting at $x = 17$ on Figure 1.2 and moving 17 steps forward we arrive at -3. For $x = 17$, equation (1.6) requires $k = 4$ to yield an integer answer. Thus we have,

$$f_4(17) = \frac{(-1)^{4-1} \cdot 3}{2^4}(17) + \frac{-2 \cdot (-1)^{4-1} + 2^4}{6 \cdot 2^4} = -\frac{3}{16}(17) + \frac{3}{16} = -3$$

2. Starting at $x = -10$ on Figure 1.2 and moving 10 steps backward we arrive at 4. For $x = -10$, equation(1.7) requires that $k = 3$ to produce an integer answer. This gives,

$$f_{-3}(-10) = \frac{(-1)^3 \cdot 3}{2^3}(-10) + \frac{4 \cdot (-1)^{3-1} + 2^3}{6 \cdot 2^3} = \frac{-3(-10)}{8} + \frac{1}{4} = 4$$

Table 1.1 below shows the resulting equations for $1 \leq k \leq 13$ and the first few values in the domain and range for equations (1.6) and (1.7).

| For Positive $x$ | Doman and Range | For Negative $x$ & 0 | Domain and Range |
|---|---|---|---|
| $f_1(x) = \frac{3x}{2}$ | In:2,4,6,8,10,12,14<br>Out:3,6,9,12,15,18,21 | $f_{-1}(x) = \frac{-3x}{2} + \frac{1}{2}$ | In:-1,-3,-5,-7,-9,-11,-13,-15<br>Out:2,5,8,11,14,17,20,23 |
| $f_2(x) = \frac{-3x}{4} + \frac{1}{4}$ | In: 3, 7, 11, 15, 19, 23<br>Out: -2, -5, -8, -11, -14, -17 | $f_{-2}(x) = \frac{3x}{4}$ | In:0, -4,-8,-12,-16,-20,-24<br>Out:0, -3,-6,-9,-12,-15,-18 |
| $f_3(x) = \frac{3x}{8} + \frac{1}{8}$ | In:5, 13, 21,29,37,45<br>Out: 2, 5, 8, 11, 14, 17 | $f_{-3}(x) = \frac{-3x}{8} + \frac{2}{8}$ | In:-2,-10,-18,-26,-34,-42<br>Out:1,4,7,10,13,16 |
| $f_4(x) = \frac{-3x}{16} + \frac{3}{16}$ | In:1,17,33,49,65,81<br>Out:0,-3,-6,-9,-12,-15 | $f_{-4}(x) = \frac{3x}{16} + \frac{2}{16}$ | In:-6,-22,-38,-54,-70,-86<br>Out:-1,-4,-7,-10,-13,-16 |
| $f_5(x) = \frac{3x}{32} + \frac{5}{32}$ | In:9,41,73,105,137, 169<br>Out:1,4,7,10,13,16 | $f_{-5}(x) = \frac{-3x}{32} + \frac{6}{32}$ | IN:-30,-62,-94,-126,-158<br>Out:3,6,9,12,15 |
| $f_6(x) = \frac{-3x}{64} + \frac{11}{64}$ | In:25,89,153,217,281, 345<br>Out:-1,-4,-7,-10,-13,-16 | $f_{-6}(x) = \frac{3x}{64} + \frac{10}{64}$ | In:-46,-110,-174,-238,-302<br>Out:-2,-5,-8,-11,-14 |
| $f_7(x) = \frac{3x}{128} + \frac{21}{128}$ | In:121,249,377,505,633<br>Out:3,6,9,12,15 | $f_{-7}(x) = \frac{-3x}{128} + \frac{22}{128}$ | In:-78,-206,-334,-462<br>Out:2,5,8,11 |
| $f_8(x) = \frac{-3x}{256} + \frac{43}{256}$ | In:185,441,697,953, 1209<br>Out: -2,-5,-8,-11, -14 | $f_{-8}(x) = \frac{3x}{256} + \frac{42}{256}$ | In:-14,-270,-526,-782,<br>Out: 0,-3,-6,-9 |
| $f_9(x) = \frac{3x}{2^9} + \frac{85}{2^9}$ | In: 313, 825, 1337, 1849<br>Out: 2, 5, 8, 11 | $f_{-9}(x) = \frac{-3x}{2^9} + \frac{86}{2^9}$ | In: -142, -654, -1166,<br>Out:1,4,7 |
| $f_{10}(x) = \frac{-3x}{2^{10}} + \frac{171}{2^{10}}$ | In: 57, 1081, 2105, 3129<br>Out: 0,-3, -6, -9 | $f_{-10}(x) = \frac{3x}{2^{10}} + \frac{170}{2^{10}}$ | In: -398, -1422, -2446,<br>Out:-1,-4,-7 |
| $f_{11}(x) = \frac{3x}{2^{11}} + \frac{341}{2^{11}}$ | In: 569, 2617, 4665, 6713<br>Out: 1, 4, 7, 10 | $f_{-11}(x) = \frac{-3x}{2^{11}} + \frac{342}{2^{11}}$ | In: -1934, -3982, -6030,<br>Out:3,6,9 |
| $f_{12}(x) = \frac{-3x}{2^{12}} + \frac{683}{2^{12}}$ | In: 1593, 5689, 9785<br>Out: -1, -4, -7 | $f_{-12}(x) = \frac{3x}{2^{12}} + \frac{682}{2^{12}}$ | In: -2958, -7054, -11150,<br>Out:-2,-5,-8 |
| $f_{13}(x) = \frac{3x}{2^{13}} + \frac{1365}{2^{13}}$ | In: 7737, 15929, 24121<br>Out: 3, 6, 9 | $f_{-13}(x) = \frac{-3x}{2^{13}} + \frac{1366}{2^{13}}$ | In:-5006, -13198, -21390<br>Out: 2, 5, 8 |

Table 1.1. Examples of input/output elements according to Lemma 1.2.

**Remark 1.2** There are many interesting observations that can be derived from Table 1.1; one of which is the repetition of the output values. Though there are an infinite number of Diphantine equations, there are only six sets of outputs. These are

i. {3,6,9,12,15,18,21,...}
ii. {-2, -5, -8, -11, -14,...}
iii. {2, 5, 8, 11, 14, 17,...}
iv. {0,-3,-6,-9,-12,-15,...}
v. {1,4,7,10,13,16,19,...}
vi. {-1,-4,-7,-10,-13,-19,...}.

**Proof of Equation (1.7):** We will prove equation (1.7); the proof of equation (1.6) is similar. Consider the integers of the form $6x + 1$ for $x \geq 0$ and $x \in \mathbb{N}$. The first iterate, $T^1$, of these numbers under the Collatz function (1.1) will all be even numbers (See Table 1.2 below).

| x | 6x+1 | $T^1$ | $T^2$ | $T^3$ | $T^4$ | $T^5$ | $T^6$ | $T^7$ | $T^8$ | $T^9$ |
|---|---|---|---|---|---|---|---|---|---|---|
| 1 | 7 | 22 | 11 | 34 | 17 | 52 | 26 | 13 | 40 | 20 |
| 2 | 13 | 40 | 20 | 10 | 5 | 16 | 8 | 4 | 2 | 1 |
| 3 | 19 | 58 | 29 | 88 | 44 | 22 | 11 | 34 | 17 | 52 |
| 4 | 25 | 76 | 38 | 19 | 58 | 29 | 88 | 44 | 22 | 11 |
| 5 | 31 | 94 | 47 | 142 | 71 | 214 | 107 | 322 | 161 | 484 |
| 6 | 37 | 112 | 56 | 28 | 14 | 7 | 22 | 11 | 34 | 17 |
| 7 | 43 | 130 | 65 | 196 | 98 | 49 | 148 | 74 | 37 | 112 |
| 8 | 49 | 148 | 74 | 37 | 112 | 56 | 28 | 14 | 7 | 22 |
| 9 | 55 | 166 | 83 | 250 | 125 | 376 | 188 | 94 | 47 | 142 |
| 10 | 61 | 184 | 92 | 46 | 23 | 70 | 35 | 106 | 53 | 160 |
| 11 | 67 | 202 | 101 | 304 | 152 | 76 | 38 | 19 | 58 | 29 |
| 12 | 73 | 220 | 110 | 55 | 166 | 83 | 250 | 125 | 376 | 188 |
| 13 | 79 | 238 | 119 | 358 | 179 | 538 | 269 | 808 | 404 | 202 |
| 14 | 85 | 256 | 128 | 64 | 32 | 16 | 8 | 4 | 2 | 1 |
| 15 | 91 | 274 | 137 | 412 | 206 | 103 | 310 | 155 | 466 | 233 |
| 16 | 97 | 292 | 146 | 73 | 220 | 110 | 55 | 166 | 83 | 250 |
| 17 | 103 | 310 | 155 | 466 | 233 | 700 | 350 | 175 | 526 | 263 |
| 18 | 109 | 328 | 164 | 82 | 41 | 124 | 62 | 31 | 94 | 47 |
| 19 | 115 | 346 | 173 | 520 | 260 | 130 | 65 | 196 | 98 | 49 |
| 20 | 121 | 364 | 182 | 91 | 274 | 137 | 412 | 206 | 103 | 310 |
| 21 | 127 | 382 | 191 | 574 | 287 | 862 | 431 | 1294 | 647 | 1942 |
| 22 | 133 | 400 | 200 | 100 | 50 | 25 | 76 | 38 | 19 | 58 |

Table 1.2

To get to the next odd number (highlighted), these $T^1$ iterates need to be divided by 2, or 4, or 8, or 16, or $2^k$ for $k \in \mathbb{N}$. We can quickly see from Table 1.2 the alternating pattern in the columns from numbers of the form $6x - 1$ when $x$ is odd to numbers of the form $6x + 1$ when $x$ is even. Also, when $x$ is odd, $T^1$ requires division by 2 to get to the next odd number which falls under $T^2$. When $x = 4n$ for n $\in \mathbb{N}$ , $T^1$ requires division by 4 to get to the next odd number which falls under $T^3$. When $x = 8n + 2$, $T^1$ requires division by 8 to get to the next odd number which falls under $T^4$, and so on.

To get directly from $x$ to the next odd number in the trajectory we multiply $x$ by 6 and add 1 to bring it back to the form $6x + 1$, then multiply $6x + 1$ by 3 and add 1 according to the Collatz function (1.1), and finally divide by $2^k$. This gives $\frac{18x+4}{2^k}$.

But $\frac{18x+4}{2^k}$ is in the form $6x - 1$ when $k$ is odd and in the form $6x + 1$ when k is even (see Table 1.2). Odd values of $k$ occur at even iterations, that is, $T^2, T^4, T^6$ and so on, while even values of $k$ occur at odd iterations.

Considering the odd values of $k$ we can add 1 and divide by 6 to get closer to the output of $x$ on Figure 1.2, which yields

$$\frac{\frac{18x+4}{2^k}+1}{6}=\frac{3}{2^k}x+\frac{2^k+4}{6\cdot 2^k}.$$

To account for the fact that $x$ is negative on the spiral for numbers of the form $(6x+1)$ we must change the sign of the coefficient of $x$ from positive to negative. This gives the final result for odd values of $k$ to be

$$-\frac{3}{2^k}x+\frac{2^k+4}{6\cdot 2^k}. \tag{1.8}$$

For even values of $k$ we subtract 1 and divide by 6 to get closer to the output of $x$ on Figure 1.2. This yields

$$\frac{\frac{18x+4}{2^k}-1}{6}$$

$$=\frac{3x}{2^k}+\frac{4-2^k}{6\cdot 2^k}.$$

Again, to account for the fact that $x$ is negative on the spiral for numbers of the form $(6x+1)$ we must change the sign of the coefficient of $x$ from positive to negative. But the outputs for even values of $k$ are of the form $6x-1$, which are positive on the spiral. To achieve a positive output, we multiply $-\frac{3}{2^k}x+\frac{4-2^k}{6\cdot 2^k}$ by $-1$ to get the final result for even values of $k$ as

$$\frac{3}{2^k}x+\frac{2^k-4}{6\cdot 2^k}. \tag{1.9}$$

If we combine equations (1.8) and (1.9) into one equation, we get the final output as $\frac{(-1)^k\cdot 3}{2^k}x+\frac{4\cdot(-1)^{k-1}+2^k}{6\cdot 2^k}$ which is equation (1.7).■

Observation. Consider equation(1.6). The constant term of each output equation is a rational number with denominator $2^k$ and numerator

$c_k=\frac{-2\cdot(-1)^{k-1}+2^k}{6}$.

If we list the values of $c_k$ we get Table 1.3.

| $k$ | 1 | 2 | 3 | 4 | 5 | 6 | 7 | 8 | 9 | 10 | 11 | 12 | 13 | … |
|---|---|---|---|---|---|---|---|---|---|---|---|---|---|---|
| $c_k = \frac{-2\cdot(-1)^{k-1}+2^k}{6}$ | 0 | 1 | 1 | 3 | 5 | 11 | 21 | 43 | 85 | 171 | 341 | 683 | 1365 | … |

Table 1.3. Numerators of the constant terms in Equation (1.6) for given k-values.

A careful look at the values of $c_k$ reveals that **if $k$ is prime and $k \geq 5$, then k divides $c_k$**. For example, 5|5, 7|21, 11|241, and 13|1365. This is not a coincidence as the sequence $c_k = \frac{-2\cdot(-1)^{k-1}+2^k}{6}$ is the Jacobsthal sequence (A001045) which has a direct link to Mersenne primes and Fermat's Little Theorem.

This intricate yet clear-cut link between the Collatz conjecture and primes is a pleasant surprise.

## 5. Graphs of the Diophantine Equations from Table 1.1

If we graph equations (1.6) and (1.7) for $k$ over the interval $[1, \infty)$ three important observations are immediately apparent:

1. the lines intersect at two common points. For equation(1.6) the common point is $\left(\frac{1}{9}, \frac{1}{6}\right)$ and for equation(1.7) the common point is $\left(\frac{2}{9}, \frac{1}{6}\right)$.
2. the lines alternate between a positive and a negative slope/coefficient;.
3. the lines have slopes/coefficients that appraoch 0 as $k$ approaches infinity.

Much could be said here. Suffice to say, if both sets of lines are graphed on the same cartesian plane, it can be observed that the lines are symmetric about $x = \frac{1}{6}$.

Figure 1.3 below shows the graphs for equation (1.6) for $x > 0$ and equation (1.7) for $x \leq 0$. Recall that equation (1.6) only accepts positive integers, while equation (1.7) accepts only negative integers and zero.

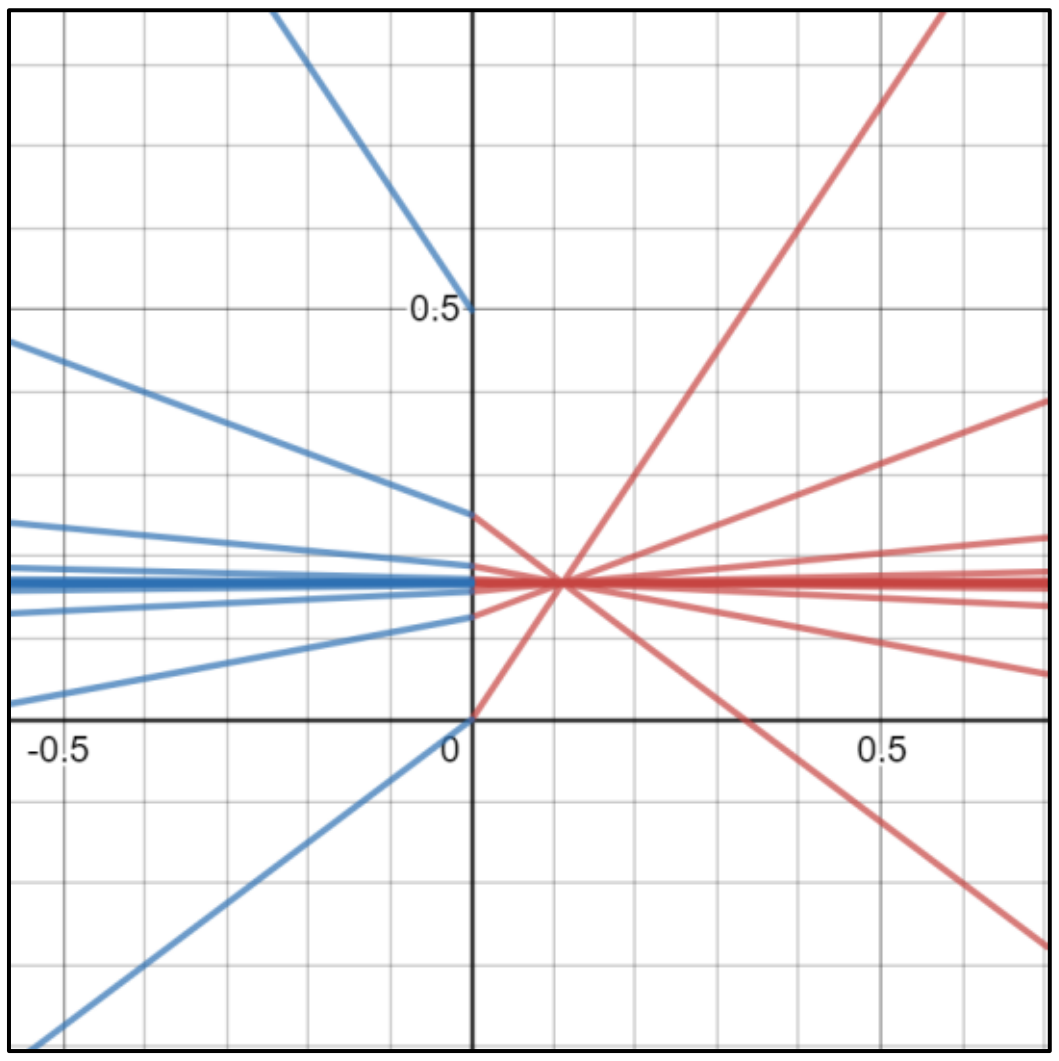

Figure 1.3. Graphs of (1.6) (red) and (1.7) (blue)

With a visual pattern to work with, a new path has opened to pinning down the Collatz Problem.

## 6. Does the Spiral Form a Connected Tree?

As stated before, a careful examination of the infinite Diophantine equations, twenty-six of which are presented in Table 1.1, will reveal that there are only six unique set of integers in the outputs. These are

i. {3,6,9,12,15,18,21,...}
ii. {-2, -5, -8, -11, -14,...}
iii. {2, 5, 8, 11, 14, 17,...}
iv. {0,-3,-6,-9,-12,-15,...}
v. {1,4,7,10,13,16,19,...}
vi. {-1,-4,-7,-10,-13,...}

As can be seen, the sets have no intersection but their union is the set of integers creating an infinitely-many-to-one mapping. That is, each value in the range is the output of infinitely many values in the domain. For example, 0 in the range is the output for the lowest value in the domain of $f_4$, $f_{10}$, $f_{16}$, ...$f_{6n+4}$ and $f_{-2}$, $f_{-8}$, $f_{-14}$, ... $f_{-6n-2}$ for $n \epsilon \mathbb{N}$. All the values that map to 0 are given by the expressions $\frac{2^{6n-2}+2}{18}$ and $\frac{4-2^{6n-4}}{18}$ for $n \epsilon \mathbb{N}$. In a similar way, all the values that map to any integer in the range can be defined by two expressions, one for positive integers and the other for negative integers and zero. It should be pointed out that the input values of all the Diophantine equations form the set of integers as well, with each integer being the input of one and only one equation. Figure 1.4 below shows the domain in the form of algebraic expressions for the range {-3, -2, -1, 0, 1, 2, 3}. Each expression generates infinite domain values for the same value in the range.

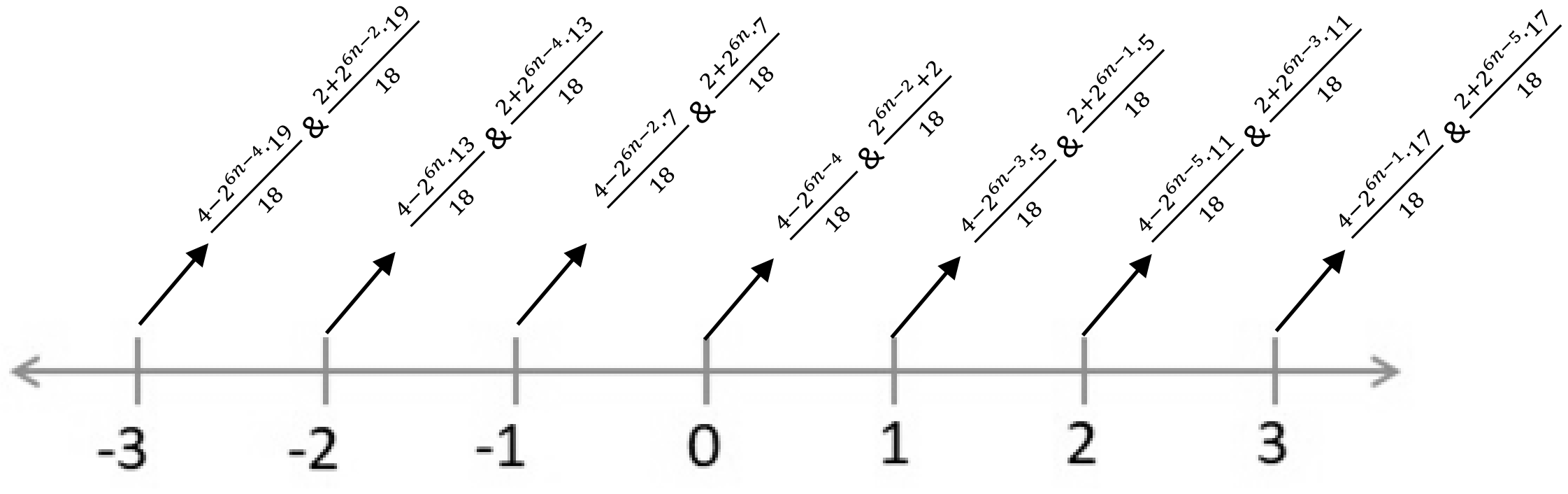

Figure 1.4. The range {-3, -2, -1, 0, 1, 2, 3} with infinite domains in the form of algebraic expressions.

Figure 1.5 below shows four out of the infinite domain values for each element in the range {-3, -2, -1, 0, 1, 2, 3}.

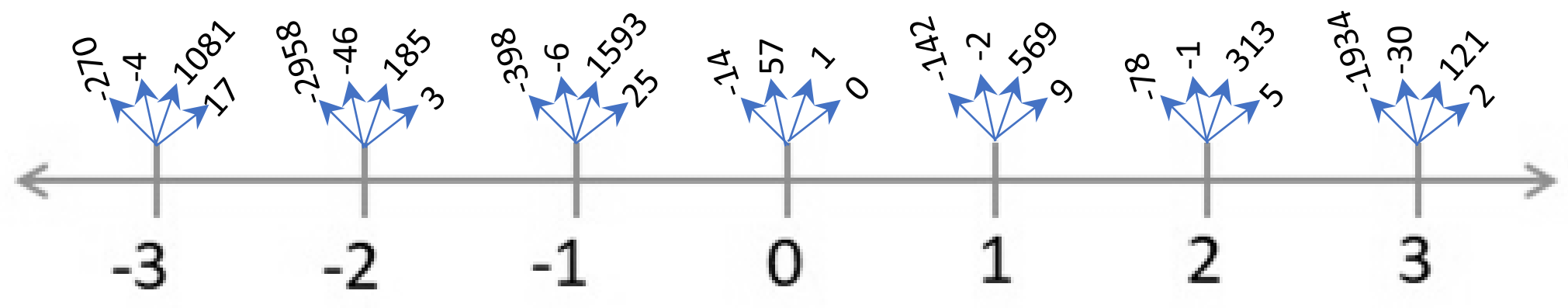

Figure 1.5. Four domain values for each element in the range {-3, -2, -1, 0, 1, 2, 3}.

If we move each value in the range to connect it with the same value in the domain, we begin to see a tree forming. Figure 1.6 below shows the tree beginning to take shape as we move just the value 2 in the range to connect it to the value 2 in the domain of another function.

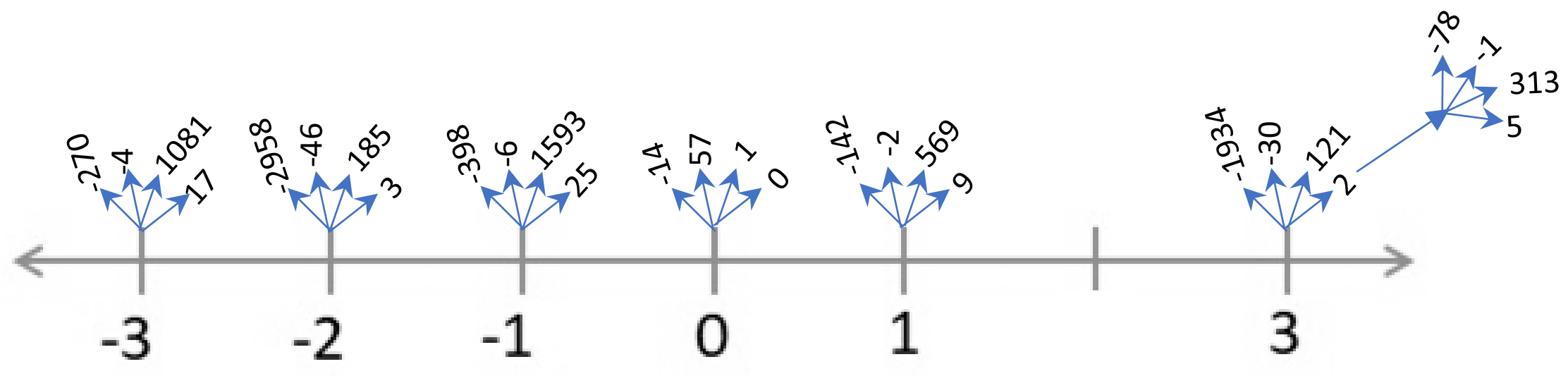

Figure 1.6. Tree shaped from connecting 2 in a range to 2 in the domain of another function.

This is function composition in which the range of one function is the domain of another. Since the domain is the set of integers and the range is the set of integers we are guaranteed to have a single value in the domain to match each value in the range as we build the tree, except for 0 which is the root and which loops back to itself.

Figure 1.7 below shows the final tree for the range {-3, -2, -1, 0, 1, 2, 3}.

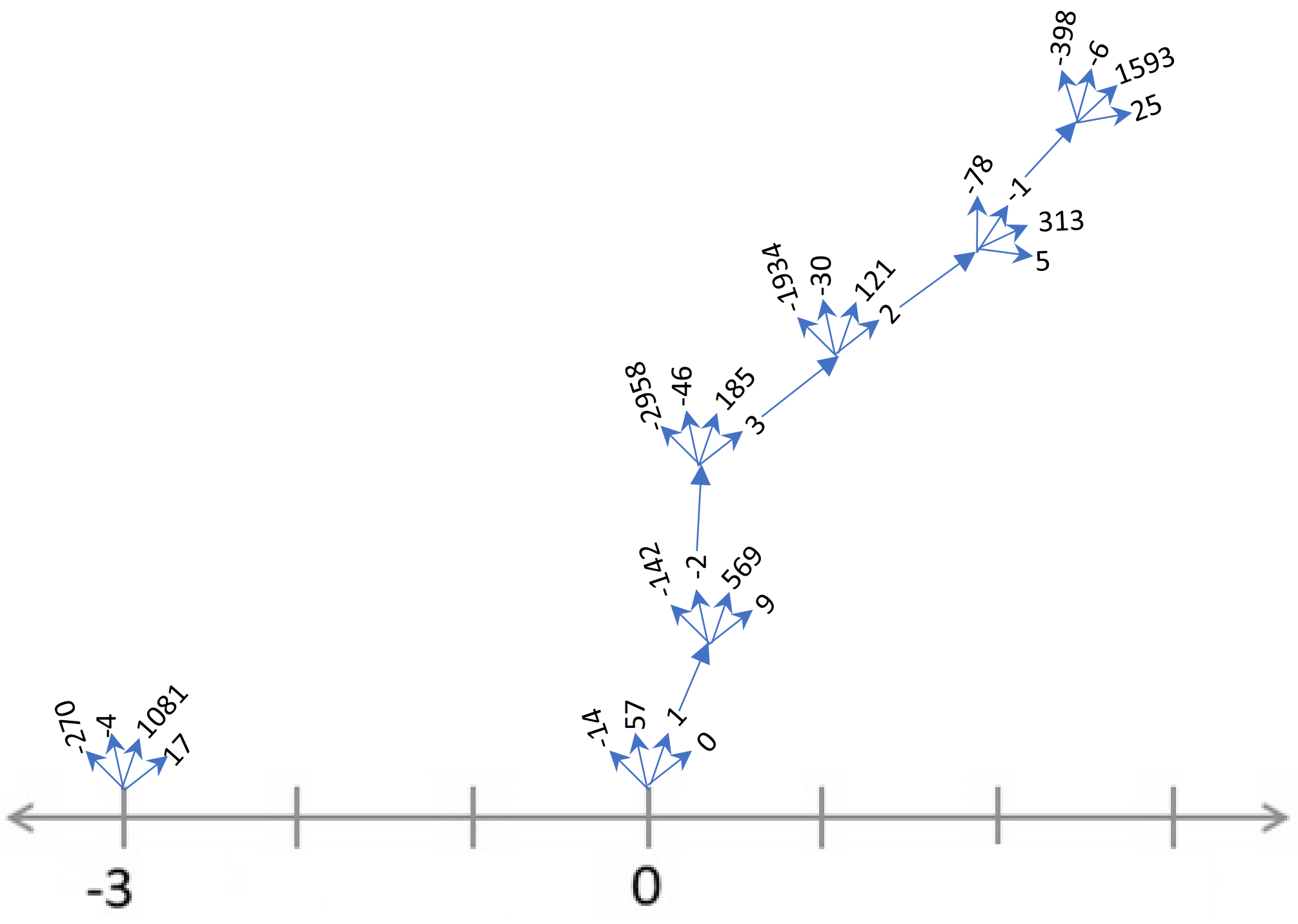

Figure 1.7. Final tree for the range {-3, -2, -1, 0, 1, 2, 3}.

We make no attempt to prove that the spiral forms a connected tree, only that there are no non-trivial cycles, which we will prove in the next section. However, we wish to point out that we can easily create similar trees from Diophantine equations that appear to contain all $\mathbb{Z}$. This means that a general method is needed to prove that all such Collatz-like trees contain all $\mathbb{Z}$.

We will now show that the spiral has no non-trivial cycles and by extension the Collatz problem.

## 7. Returning to a Previous Number on the Spiral

From Table 1.1 it can be seen that each function on the right-hand side of the table will immediately return its input values to a previous number on the spiral. This is because the inputs are all negative values which move backward on the spiral. Here we wish to emphasize that the meaning of "previous number on the spiral" is as it reads; any number that comes before the number in question on the spiral. Specifically, the phrase "previous number on the spiral" does not always mean a number with a lower value – it means literally "preceding the number in question on the spiral". There are numbers with greater value when converted back to the form $6x \pm 1$ that precede smaller values of the form $6x \pm 1$ on the spiral. For example, 11 which is represented by 2 on the spiral comes before 7 which is represented by -1 on the spiral, even though 11 is larger than 7.

Now let us return to Table 1.1. On the left-hand side of the table, each function below $f_2$ also returns its input values to previous numbers on the spiral. This is somewhat intuitive as division by 8, 16, 32 and so on represents falling at least one level on the spiral. The only two functions that do not immediately produce previous values on the spiral are $f_1$ and $f_2$. In fact, for the trajectory of any number to climb to higher values it must go through $f_1$ and/or $f_2$.

**Proposition 1.2** No non-trivial cycle can exist without an element from $f_1$ and/or $f_2$.

We therefore now shift our focus to $f_1$ and $f_2$.

Starting at either $f_1$ or $f_2$, if the trajectory of a number goes to a previous number on the spiral, then an infinite arithmetic sequence of numbers takes the same path to previous numbers on the spiral. For example, consider the composition that starts at $f_2$ then moves to $f_{-2}$ and then to $f_{-1}$:

$$f_{-1}f_{-2}f_2(x) = -\frac{3}{2}\left(\frac{3}{4}\left(-\frac{3}{4}x + \frac{1}{4}\right) + 0\right) + \frac{1}{2} = \frac{27}{32}x + \frac{7}{32}.$$

Note that the coefficient of $x$ is positive and less than 1. The smallest number in $f_2$ that satisfies this composition is 27. Now $f_{-1}f_{-2}f_2(27) = \frac{27}{32}(27) + \frac{7}{32} = 23$. Hence the composition returns the input value of 27 to a previous number on the spiral, 23. What is interesting is that all numbers of the form $23 + 32\mathrm{k}$ for $k \in \mathbb{N}$ also belong to the domain of $f_2$ and backtrack to previous numbers on the spiral following the same path and ending at the same function under $f_{-1}f_{-2}f_2(x) = \frac{27}{32}x + \frac{7}{32}$. For example,

$f_{-1}f_{-2}f_2(27 + 32) = f_{-1}f_{-2}f_2(59) = \frac{27}{32}(59) + \frac{7}{32} = 50$. That is, 59 from $f_2$ goes to 50 in $f_{-1}$, a previous number on the spiral. Note that the common difference in $23 + 32\mathrm{k}$ is simply the denominator of the composite function $f_{-1}f_{-2}f_2(x)$, and the common difference in the output values is the numerator of the coefficient. This leads to Proposition 1.3.

**Proposition 1.3** Any composition starting at $f_1$ or $f_2$that has a positive coefficient that is less than one will return its input values to previous numbers on the spiral.

On the other hand, if the coefficient is negative, then it must be greater than $-\frac{1}{2}$ to return its input values to previous numbers on the spiral. This is because the absolute value of any negative number on the spiral is only half the positive number that comes just before. Take for example the composition $f_4f_1f_1f_{-1}f_2(x) = -\frac{243}{512}x + \frac{69}{512}$, which brings the inputs 211, 723, 1235… to -100, -343, -586, … respectively. Note that the absolute value of each negative output is less than half the corresponding input. This means that the trajectory of each input has reached a previous number on the spiral. Dipping lower in the table, for example instead of $f_4f_1f_1f_{-1}f_2(x)$ we have $f_6f_1f_1f_{-1}f_2(x)$, which also guarantees that each input will return to a previous value on the spiral. This motivates the following definition.

**Definition 1.1** A Sink Composition is any composition that produces previous values on the spiral.

Sink Compositions beg the question: What powers of 3 and 2 will give a positive coefficient that is less than one, or a negative coefficient that is greater than $-\frac{1}{2}$? The solutions to $\frac{3^p}{2^q} < 1$ (or $p < \frac{qlog(2)}{\log(3)}$ when rearranged) and $\frac{-3^p}{2^q} > -\frac{1}{2}$ (or $p < \frac{xlog(2)-\log(2)}{\log(3)}$ when rearranged) answer the question. The graphs of both inequalities are shown below in Figure 1.8 without the shadings but with the understanding that the solutions are all the lattice points that fall below the lines.

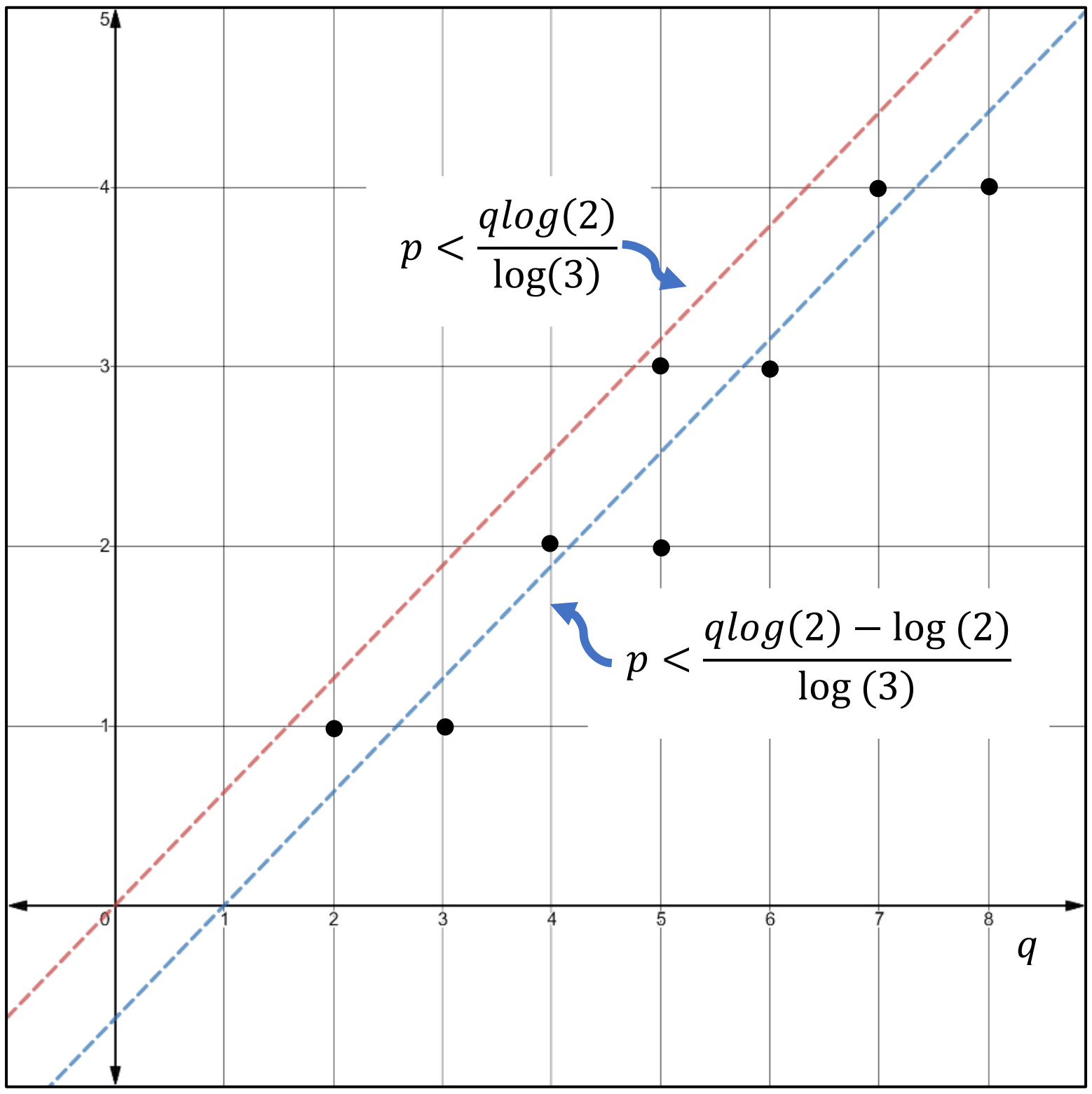

Figure 1.8. Graph of $p < \frac{qlog(2)}{\log(3)}$ and $p < \frac{xlog(2)-\log(2)}{\log(3)}$.

Since the p-values are the powers of 3 which represent the number of steps in the spiral trajectory, the graph of $p < \frac{qlog(2)}{\log(3)}$ reveals that for positive coefficients to be less than 1, if there is only one step then the minimum allowable power of 2 is 2. That is, if the numerator is $3^1$ then the denominator must be $2^2$ or higher for the fraction $\frac{3^p}{2^q}$ to be less than 1. This is represented by the lattice point (2, 1) on the graph. If there are two steps, then the minimum allowable power of 2 is 4, represented by (4, 2); if there are three steps, then the minimum allowable power of 2 is 5, represented by (5,3); and so on.

Likewise, the graph of $p < \frac{qlog(2)-\log(2)}{\log(3)}$ reveals that for negative coefficients to be greater than $-\frac{1}{2}$, if there is only one step then the minimum allowable power of 2 is 3. That is, if the numerator is $-3^1$ then the denominator must be $2^3$ or higher for the fraction $\frac{-3^p}{2^q}$ to be greater than $-\frac{1}{2}$. This is represented by the lattice point (3, 1). If there are two steps, then the minimum allowable power of 2 is 5, represented by (5, 2); if there are three steps, then the minimum allowable power of 2 is 6, represented by (6, 3); and so on.

## 8. What Composition can Output its Input to Complete a Cycle?

Recall from Proposition 1.3 that no non-trivial cycle can exist without an element from $f_1$ or $f_2$. With that in mind, we could classify compositions of the Diophantine equations/functions of the form given in Lemma 1.2 originating from either $f_1$ or $f_2$ into four categories:

1. Category 1: Compositions with a positive coefficient of $x$ and a positive constant,
2. Category 2: Compositions with a positive coefficient of $x$ and a negative constant,
3. Category 3: Compositions with a negative coefficient of $x$ and a positive constant, and
4. Category 4: Compositions with a negative coefficient of $x$ and a negative constant.

Note that a composition that results in a cycle will produce an output that is equal to the input for one of the feasible input values. That is, $\frac{\pm 3^m}{2^n}x_i \pm \frac{c}{2^n} = x_i$, for some input value $x_i$.

Starting from either $f_1$ or $f_2$ where the inputs are positive, it is easy to deduce that Category 4 compositions will produce negative output values. Such compositions could not result in a cycle from $f_1$ or $f_2$, since the outputs are negative, but the inputs are positive.

For compositions in Category 1, if the coefficient is greater than 1, then the output would be larger than the input and a cycle would not be possible. If the coefficient is less than 1, then it is a sink composition and returns the input to a previous number on the spiral.

For Category 3 compositions where the coefficient is negative and the constant is positive, the outputs will always be negative. This is because (i) there is no Diophantine equation in Table 1.1 with a range containing a mix of both positive and negative numbers, which means the output values are either all positive or all negative and zero; (ii) no matter the size of the constant there exists a large enough value in $f_1$ or $f_2$ that will cause the function to produce a negative output since if $x_i$ satisfies the Diophantine composition $\frac{-3^m}{2^n}x + \frac{c}{2^n}$ then $x_i + 2^n \cdot k$ also satisfies the composition ($c, k, m$ *and* $n$ are whole numbers). For suppose there exists $c$ large enough and $x_i$ is small enough so that $\frac{-3^m}{2^n}x_i + \frac{c}{2^n}$ produces a positive output, but at the same time $x_i + 2^n \cdot k$ can be made large enough so that the composition produces a negative value. This is a contradiction of the fact that none of the Diophantine equations of the form given in Lemma 1.2 produces both positive and negative outputs. Hence, we conclude that a composition with a negative coefficient and a positive constant can only produce negative outputs or zero. We therefore conclude that Category 3 composition could not return a number to itself to complete a cycle.

Having shown that compositions originating from $f_1$ or $f_2$ that fall into Categories 1, 3 and 4 cannot return an element in the domain of $f_1$ or $f_2$ to itself, we are left with only Category 2 compositions to deal with. Category 2 compositions seem like the ideal compositions for a cycle. If the coefficient $\frac{3^m}{2^n}$ is slightly greater than 1 then adding a negative constant can decrease the product $\frac{3^m}{2^n}x_i$, which would be slightly greater than $x_i$, to equal to $x_i$ thus completing a cycle; that is, $\frac{3^m}{2^n}x_i - \frac{c}{2^n} = x_i$. But the important question here is: Can such a composition be produced? If no composition exists that falls into Category 2, then a cycle cannot be created containing an element from the domain of $f_1$ or $f_2$. This would also mean that no non-trivial cycles exist on the spiral; we address this question in the next section. Recall that by Proposition 1.3 all other functions sink immediately and therefore without $f_1$ and/or $f_2$ the spiral trajectory would converge to zero.

## 9. Do Category 2 Compositions with Positive Coefficients and Negative Constants Exist?

Half the Diophantine equations in Table 1.1 have positive x-coefficients and positive or zero constants. These functions have the attributes of Category 1 compositions. The other half of the Diophantine equations have a negative x-coefficient and a positive constant which are attributes of Category 3 Compositions. We can create Category 4 compositions with negative coefficients and negative constants by repeatedly going back and forth between $f_1$ and $f_2$. For example, the composition $f_2f_{-1}f_2f_{-1}f_2f_{-1}f_2(x) = -\frac{3^7}{2^{11}}x - \frac{139}{2^{11}}$. We could also bring the composition to $f_1$ to increase the constant faster, then bring the composition back through $f_2$ to create a Category 4 composition. For example,

$$f_2f_1f_1f_1f_1f_1f_1f_1f_1f_{-1}f_2(x) = -\frac{3^{11}}{2^{13}}x - \frac{17635}{2^{13}} = -\frac{177147}{8192}x - \frac{171635}{8192}$$

which has domain $2503 + 2^{13} \cdot k$ and range $-177147k + 123019$.

If a composition with a positive coefficient and a negative constant has a coefficient that is less than 1, then it is a sink composition and produces previous numbers on the spiral in comparison to the input values; hence it cannot create a cycle. This means that if a composition with a positive x-coefficient and a negative constant were to output its input the coefficient must be larger than 1.

Let us call the desired Category 2 equation with a positive x-coefficient and a negative constant Equation 1, that is,

$$\frac{3^{m+1}}{2^{n_1+n_2}}x_i \overbrace{- \frac{3}{2^{n_2}} \cdot \frac{c_1}{2^{n_1}} + \frac{c_2}{2^{n_2}}}^{\text{Negative}} = x_i \quad \text{...... Equation 1}$$

where the constant is negative. To achieve Equation 1, it must be derived from a Category 3 equation of the form

$$g(x) = -\frac{3^m}{2^{n_1}}x_i + \frac{c_1}{2^{n_1}} \text{ .... Equation 2.}$$

To see this, note that when a Category 4 composition is substituted into any of the Diophantine functions, the result is either a Category 1 composition, a Category 3 composition or a Category 4. Also, when a Category 1 composition is substituted into any of the Diophantine functions, the result is a Category 1 composition, a Category 3, or a Category 4 composition.

Specifically, if we take Equation 2 (which produce negative outputs) and substitute it into a Category 3 function with a negative coefficient and a positive constant from the right-hand side of Table 1.1, say $f_{-n_2}(x) = -\frac{3}{2^{n_2}}x + \frac{c_2}{2^{n_2}}$, then we get Equation 1 as follows:

$$f_{-n_2}g(x) = -\frac{3}{2^{n_2}}\left(-\frac{3^m}{2^{n_1}}x_i + \frac{c_1}{2^{n_1}}\right) + \frac{c_2}{2^{n_2}} = \frac{3^{m+1}}{2^{n_1+n_2}}x_i \overbrace{- \frac{3}{2^{n_2}} \cdot \frac{c_1}{2^{n_1}} + \frac{c_2}{2^{n_2}}}^{\text{Negative}}$$

For Equation 1 to be derived from Equation 2, the constant in Equation 2 must be greater than $\frac{1}{3}$. To see this, note that the constant in Equation 1 is negative, that is,

$$-\frac{3}{2^{n_2}} \cdot \frac{c_1}{2^{n_1}} + \frac{c_2}{2^{n_2}} < 0$$

which yields

$$\frac{c_1}{2^{n_1}} > \frac{c_2}{3}.$$

We now seek the smallest value for $c_2$. The smallest numerator of a constant belonging to a function with a negative x-coefficient and a positive constant on the right-hand side of Table 1.1 is 1. Hence,

$$\frac{c_1}{2^{n_1}} > \frac{1}{3}.$$

Indeed, if a function with a constant greater than $\frac{1}{3}$ is substituted into $f_{-1}$ then the resulting function will have a negative constant, which is what we desire in Equation 1. We could also get $\frac{c_1}{2^{n_1}} > \frac{1}{3}$ by substituting the constant from Equation 2 into $f_{-1}$ like so

$$-\frac{3}{2}\cdot\frac{c_1}{2^{n_1}}+\frac{1}{2}<0$$

which yields $\frac{c_1}{2^{n_1}} > \frac{1}{3}$. There is an important point here to be made which we will state as Proposition 1.4.

**Proposition 1.4** Whenever a composite function is created by substituting one Diophantine equation into another, the constant of the composite function is simply the answer you get when the constant of the substituted equation replaces the variable $x$ in the receiving equation.

For example, if we wish to find the constant for the composite function $f_{-1}f_2(x)$ where $f_{-1}(x) = -\frac{3}{2}x + \frac{1}{2}$ and $f_2(x) = -\frac{3}{4}x + \frac{1}{4}$, then we can simply replace $x$ in $f_{-1}$ with $\frac{1}{4}$ to get $-\frac{3}{2}\left(\frac{1}{4}\right) + \frac{1}{2} = \frac{1}{8}$. This point is important because the graph of each Diophantine equation are all the possible constants the equation can produce. Take the graph of $f_2$, for example, shown below in Figure 1.9.

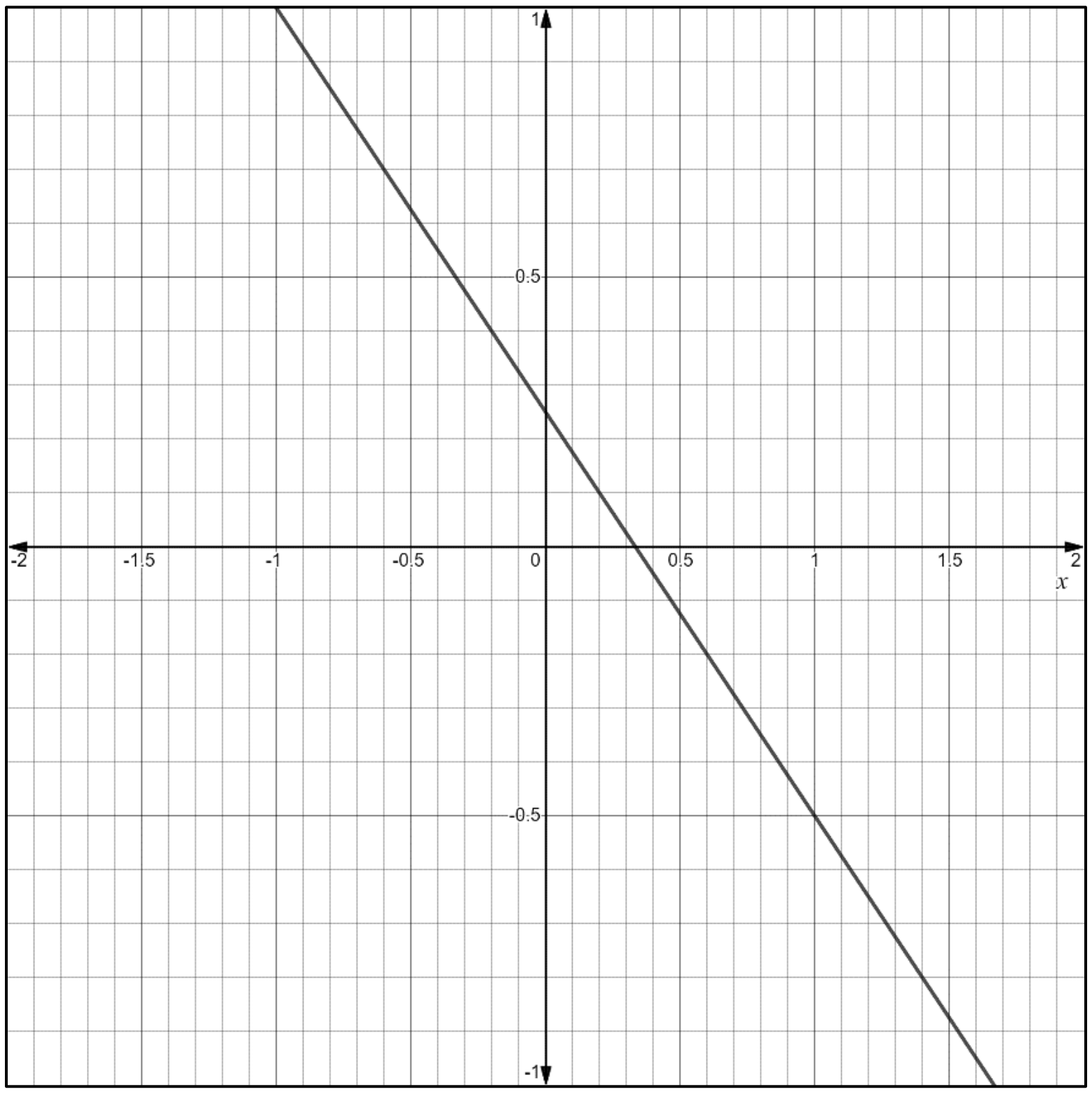

Figure 1.9. Graph of $f_2(x) = -\frac{3}{4}x + \frac{1}{4}$.

Three critical segments can be highlighted: $x < 0,\ 0 < x < \frac{1}{3}$ and $x > \frac{1}{3}$. If a function were to be substituted into $f_2$ with a constant greater than $\frac{1}{3}$, then the constant of the composite function would be negative.

i. If a function were to be substituted into $f_2$ with a constant between 0 and $\frac{1}{3}$, then the constant of the composite function would be positive but less than the constant of $f_2$.

ii. If a function were to be substituted into $f_2$ with a constant less than 0, then the constant of the composite function would be positive but greater than the constant of $f_2$.

Now the question becomes: Can we compose Equation 2 with $\frac{c_1}{2^{n_1}} > \frac{1}{3}$? We will approach this question by looking at four possibilities.

1. Can Equation 2 be composed, or outputted, by $f_{-1}, f_{-3}, f_{-5}, \ldots f_{-2k+1}$ where $k \in \mathbb{N}$? The simple answer is no. Functions with subscripts of the form $-2k+1$ all produce positive values, while Equation 2 produces negative values.
2. Can Equation 2 be composed, or outputted, by $f_{-2}, f_{-4}, f_{-6}, \ldots f_{-2k}$ where $k \in \mathbb{N}$? For functions with subscripts of the form $-2k$ to compose Equation 2, the input function must have a negative coefficient and a constant that's less than $\frac{1}{3}$. (If the constant were greater than $\frac{1}{3}$ then the input function would already be in the form of Equation 2). So, can functions with subscripts of the form $-2k$ receive a constant that is less than $\frac{1}{3}$ and output a constant that is greater than $\frac{1}{3}$? To answer this question, we turn our attention to the graphs of these functions shown in Figure 1.10 below.

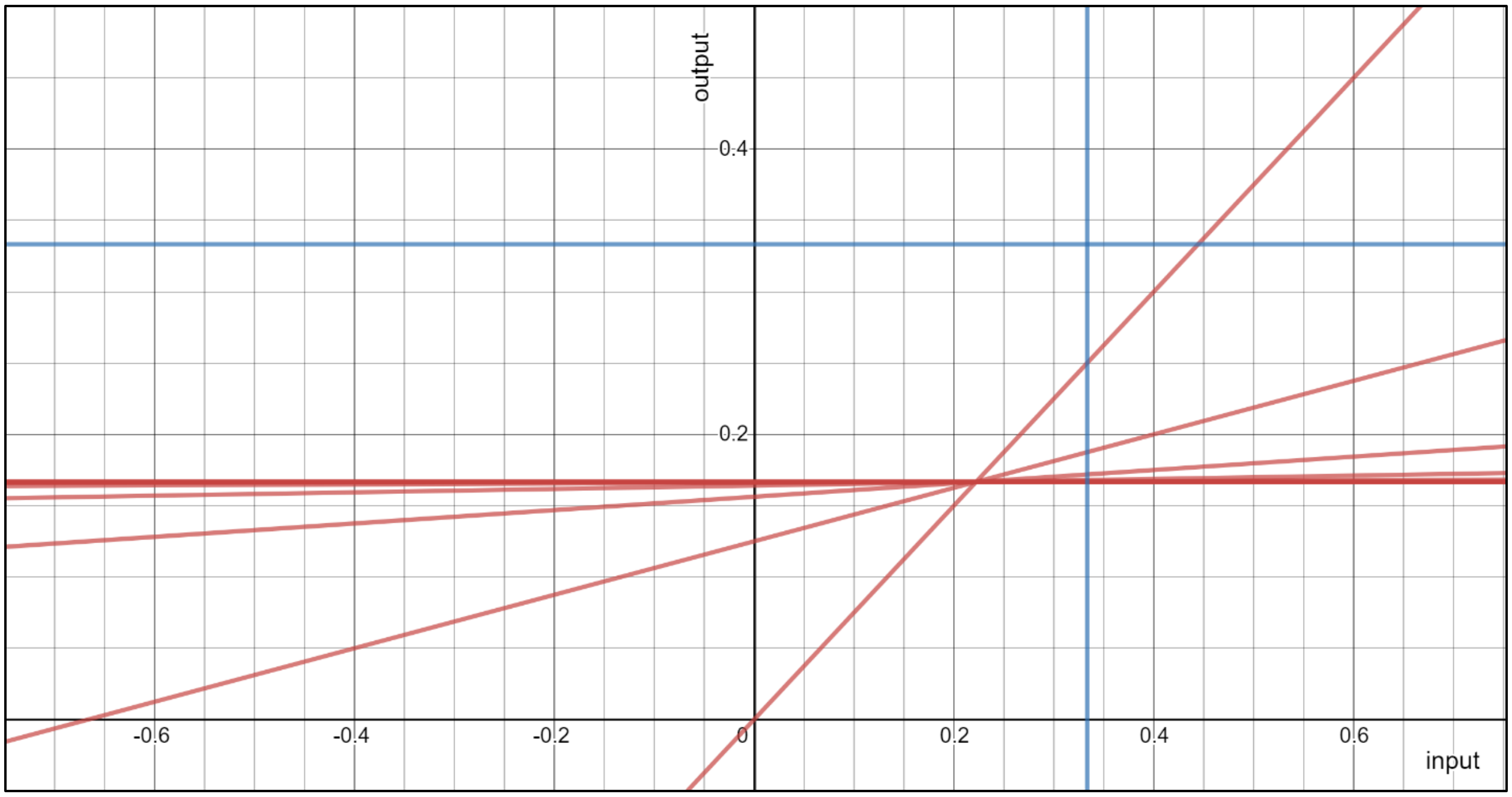

Figure 1.10. Graphs of $f_{-2k}$ where $k \in \mathbb{N}$, $x = \frac{1}{3}$ and $y = \frac{1}{3}$.

Functions with subscript $-2k$ are an infinite set of functions with coefficients/slopes that approach 0 as $k \to \infty$. The lines $x = \frac{1}{3}$ and $y = \frac{1}{3}$ are also added for reference. Note that when $x < \frac{1}{3}$ all the lines have values that are also less than $\frac{1}{3}$. This means that functions with subscripts of the form $-2k$ cannot compose, or output, Equation 2.

3. Can Equation 2 be composed, or outputted, by $f_1, f_3, f_5, \ldots f_{2k-1}$ where $k \in \mathbb{N}$? The simple answer is no. Like functions with subscripts $-2k+1$, functions with subscripts of the form $2k-1$ all produce positive values. Since Equation 2 produces negative values it cannot be composed, or outputted, by a function with subscript $2k-1$.
4. With Equation 2 not producible by any of the previous three forms of functions, it must be producible by functions with subscripts $2k$, otherwise it is impossible to compose Equation 2 using the Diophantine equations. Functions with subscripts $2k$ have a negative coefficient and a positive constant. These functions receive functions with a positive coefficient and a positive constant. (The constant cannot be negative otherwise the input function would be in the form of Equation 1 which is what we are ultimately trying to produce). A careful look at the graphs for functions with subscripts $2k$ (see Figure 1.11 below) reveals that for any positive constant inputted into these functions, they output either a positive constant that is less than $\frac{1}{3}$ or a negative constant, neither of which we seek.

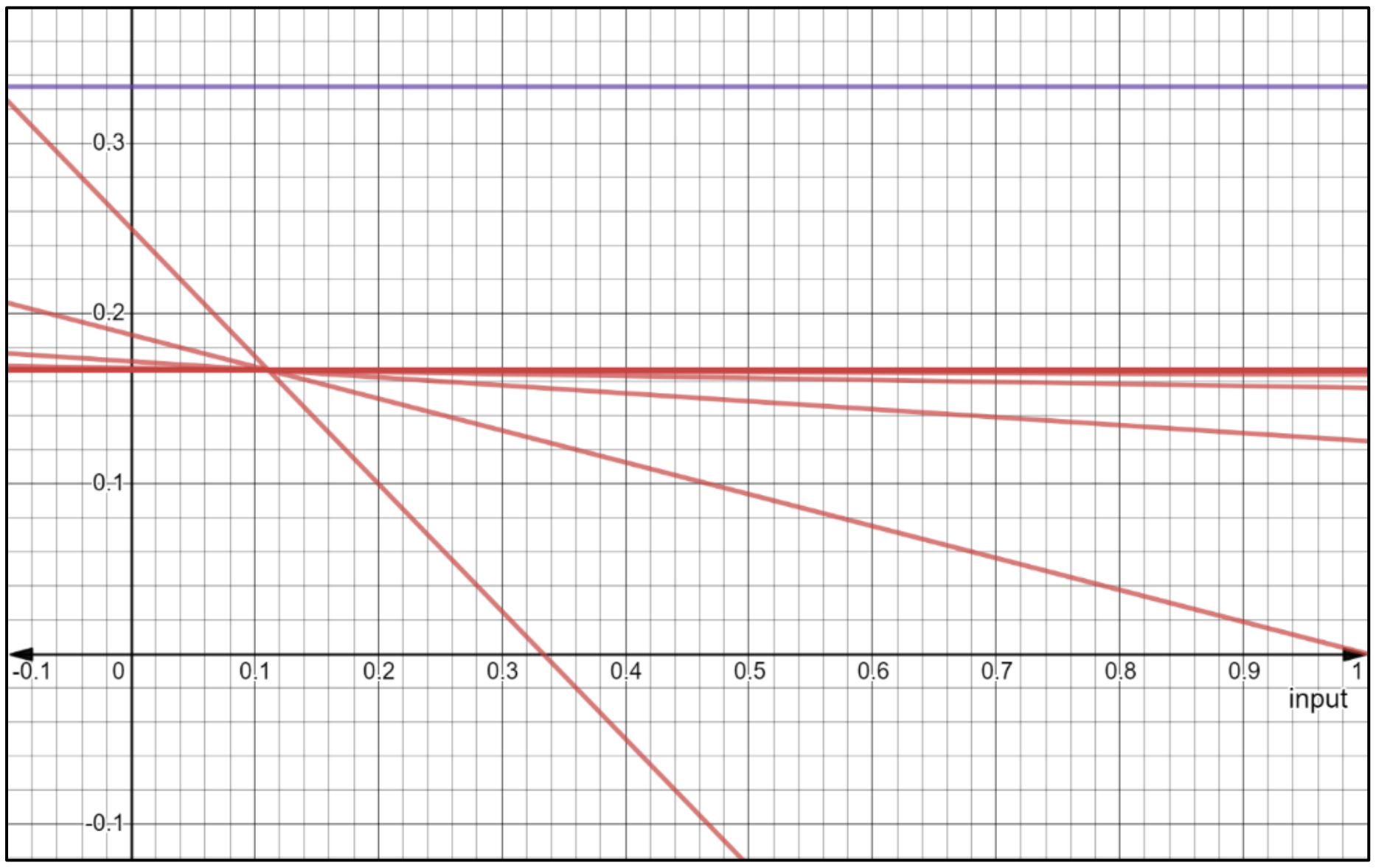

Figure 1.11. Graphs of $f_{2k}$ where $k \in \mathbb{N}$ and $y = \frac{1}{3}$.

Having examined all the Diophantine equations to see which, if any, can produce Equation 2 and finding none, we can conclude that Equation 2 is not producible, and by extension, neither is Equation 1. Thus, a cycle containing an element from $f_1$ or $f_2$ is not possible. This means that a composition resulting in a non-trivial cycle is not possible; the spiral contains no non-trivial cycles. ■

## 10. The Collatz Conjecture has no Non-Trivial Cycles

With no non-trivial cycles on the spiral and by extend in the Collatz function (1.1), we are one step closer to solving the Collatz Problem.

---

## Acknowledgements

I am grateful for the insightful comments provided by Dr. Dominic D. Clemence, professor in the department of Mathematics and Statistics at NC A&T State University. His knowledge and exacting attention to detail have saved me from many errors. I am also grateful for the feedback provided by Kells Hall, Alden Campbell, and Marvin Aidoo.